\documentclass[J,12pt]{amsart}
\usepackage{epsfig,epsf,amsmath,amssymb}
\newtheorem{theorem}{Theorem}[section]
\newtheorem{corollary}[theorem]{Corollary}
\newtheorem{definition}{Definition}

\newtheorem{proposition}[theorem]{Proposition}
\newtheorem{remark}[theorem]{Remark}

\begin{document}
\title[Remarks on pseudocovering spaces]{Pseudocovering and digital covering spaces}
\author[label1]{Sang-Eon Han}
\address[label1]{Department of Mathematics Education, Institute of Pure and Applied Mathematics\\
	Jeonbuk National University, Jeonju-City Jeonbuk, 54896, Republic of Korea\\
	e-mail address:sehan@jbnu.ac.kr, Tel: 82-63-270-4449.}

\begin{abstract} The notions of a local $(k_0,k_1)$-isomorphism and a weakly local $(k_0,k_1)$-isomorphism play crucial roles in 
	developing a digital $(k_0,k_1)$-covering space and a pseudo-$(k_0,k_1)$-covering space, respectively.
	In relation to the study of pseudo-$(k_0,k_1)$-covering spaces, since there are some works to be refined and improved in the literature, 
	the recent paper \cite{H10}
	improved and corrected some mistakes occurred in the literature.
	One of the important things is that the notion of a pseudo-$(k_0,k_1)$-covering map in \cite{H6,H9} was revised to be more broadened in \cite{H10}.
	Thus this new version is proved to be equivalent to a weakly local $(k_0,k_1)$-isomorphic surjection \cite{H10}. The present paper contains some works in \cite{H10} and
	we  only deals with $k$-connected digital images $(X, k)$.
\end{abstract}


\maketitle

\newpage

\section{\bf Introduction}\label{s1}

The notion of a pseudo-$(k_0, k_1)$-covering space was initially introduced in 2012 \cite{H6}.
Indeed, it was intended to make a digital $(k_0, k_1)$-covering space in \cite{H1,H2,H3,H5} more generalized and broader.
Hence it was defined by using three conditions among which two of them, i.e., the conditions (1) and (2) for a  pseudo-$(k_0, k_1)$-covering space (see Definition 4), are equal to those for a digital $(k_0, k_1)$-covering space(see Definition 5).
Meanwhile, the other condition (3) for a pseudo-$(k_0, k_1)$-covering space is different from the condition (3) for a digital $(k_0, k_1)$-covering space (see Definitions 4 and 5 in the present paper). 
To be specific, the former was defined by using the notion of a weakly local ($WL$-, for brevity)$(k_0, k_1)$-isomorphism and the latter was characterized by using the concept of a local  $(k_0, k_1)$-isomorphism. 
Thus these two conditions (3) are quite different from each other.
However, when combining the two conditions (1) and (2) with each of the conditions (3), a pseudo-$(k_0, k_1)$-covering space implies to a digital  $(k_0, k_1)$-covering space (see Theorem 3.3 of \cite{P1}).
Probably, in \cite{H6}, there seems to be a gap between the author's intension for establishing a pseudo-$(k_0, k_1)$-covering space and the mathematical presentation of it.
Hence the recent paper \cite{H10} revised the original version of the condition (1) for a pseudo-$(k_0, k_1)$-covering space (see Definition 4 in the present paper) to finally make a distinction between a digital $(k_0, k_1)$-covering and a revised version of the original version of a pseudo-$(k_0, k_1)$-covering map in Definition 4. 
In detail, we will shortly see the revised version of Definition 4 via Definition 4.1 of \cite{H10}.\\

The present paper only deals with $k$-connected digital images unless otherwise stated and often uses the notion $:=$ to introduce some terms. \\

The four papers \cite{H6,H9,P1,P2} studied various properties of $\lq\lq$pseudocovering spaces".
Since the two of them \cite{H6,H9} have some errors and the others \cite{P1,P2} also have some mistakes relating to the map in (1.1) below, 
the recent paper \cite{H10} corrected and improved them, which makes them so clear.
More precisely, with the original version of a  pseudo-$(k_0, k_1)$-covering map (see Definition 4 in the present paper), 
the map $p$ in (1.1) below is not a  pseudo-$(2,k)$-covering map (see Proposition 3.2 of \cite{P1}).
Since there are some errors in the proof of  Proposition 3.2 of \cite{P1}, we note that the paper \cite{H10} corrected it.
 $$ 
\left \{\aligned 
& p: ({\mathbb Z}^+, 2) \to SC_k^{n,l}:=(x_i)_{i \in [0, l-1]_{\mathbb Z}}\,\,\text{defined by}\\
& p(t)= x_{t(mod\,l)},\,\text{where}\,\,{\mathbb Z}^+:=[0, \infty)_{\mathbb Z}:=\{t\in {\mathbb Z}\,\vert \, t \geq 0\}.
\endaligned
\right\} \eqno (1.1) 
$$
Hence the recent paper \cite{H10} fully explained the process of a non-pseudo-$(2,k)$-covering map of $p$ in (1.1).

Next, we also note that there are some mistakes on the identity of (4.2) in Proposition 4.4 and Corollary 4.5 of \cite{H9}. 
The paper \cite{H10} pointed out these defects and corrected them and verified that a $WL$-$(k_0, k_1)$-surjection is not equivalent to a pseudo-$(k_0, k_1)$-covering map followed by Definition 4 in the present paper. \\

To sum up, the recent paper \cite{H10} did corrections and improvements, as follows:\\
(1) Corrections of the map $p$ of (4.1) in the proof of Remark 4.3(2) of \cite{H9}.\\
(2) Corrections of the identity of (4.2) of Proposition 4.4 and Corollary 4.5 of \cite{H9}.\\
(3) Revision of the notion of a pseudo-$(k_0, k_1)$-covering space of Definition 4.\\
(4) Correction of the proof of Proposition 3.2 of \cite{P1} and related works in \cite{P2}.\\
(5) Improvement of the proof of Theorem 3.3 of \cite{P2}. \\
In addition, we confirm that the example in (1.1) now becomes an example for the revised version of a pseudo-$(2, k)$-covering space in \cite{H10}.

\section{\bf Preliminaries}
In relation to the study of some properties of a pseudo-$(k_0, k_1)$-covering space, to make the paper self-contained, 
we will refer to some notions.
Naively, a digital image $(X, k)$ can be considered to
be a set $X\subset {\mathbb Z}^n$ with one of the $k$-adjacency of
${\mathbb Z}^n$  from (2.1) below (or a digital $k$-graph on ${\mathbb Z}^n$ \cite{H4}).
Indeed, the papers \cite{KR1,R1} considered $(X, k), X \subset {\mathbb Z}^n, n \in \{1,2,3\}$, with $2$-adjacency on ${\mathbb Z}$, $4,8$-adjacency on ${\mathbb Z}^2$, and  $6, 18, 26$-adjacency on ${\mathbb Z}^3$.
As the generalization of the low dimensional cases, the digital $k$-adjacency relations (or digital $k$-connectivity) for $X \subset {\mathbb Z}^n, n \in {\mathbb N}$, were initially established in \cite{H7} (see also \cite{H1,H2,H3}), as follows:\\
For a natural number $t$, $1 \leq t \leq n$, the distinct points
$p = (p_1, p_2, \cdots, p_n)$  and $q =(q_1, q_2, \cdots, q_n)\in {\mathbb Z}^n$
are $k(t, n)$-adjacent if at most $t$ of their coordinates  differ by $\pm1$ and the others~coincide.
Indeed, the numbers of $t$ and $n$ of $k(t, n)$ above is very important. For instance, on
${\mathbb Z}^2$, two types of digital $k$-adjacencies exist such as  $k(1,2)=4$ and $k(2,2)=8$. 
Meanwhile, on
${\mathbb Z}^4$, four kinds of digital $k$-adjacencies exist such as $k(1,4)=8$, $k(2,4)=32$, 
$k(3,4)=64$, $k(4,4)=80$.
Then, even though  the $8$-adjacency are used on both ${\mathbb Z}^2$ and ${\mathbb Z}^4$, 
using  $k(2,2)=8$ and $k(1,4)=8$, we can make a distinction between them efficiently.
According to this statement,  the well-presented $k(t, n)$-adjacency relations (or digital $k$-connectivities) of ${\mathbb Z}^n, n \in {\mathbb N}$, are formulated \cite{H7} (see also \cite{H3}) as follows:
$$k:=k(t,n)=\sum_{i=1}^{t} 2^{i}C_{i}^{n}, \text{where}\,\, C_i ^n:= {n!\over (n-i)!\ i!}. \eqno(2.1)$$
Based on the $k$-adjacency relations of ${\mathbb Z}^n$ in (2.1), $n \in {\mathbb N}$, we will call the pair $(X, k)$ a  digital image on ${\mathbb Z}^n$, $X \subset {\mathbb Z}^n$.\\

A simple closed $k$-curve (or simple $k$-cycle) with $l$ elements in ${\mathbb Z}^n, n \geq 2$, denoted by $SC_k^{n,l}$ \cite{H3,KR1}, $l(\geq 4) \in {\mathbb N}$, is defined to be the set $(x_i)_{i \in [0, l-1]_{\mathbb Z}}\subset {\mathbb Z}^n$ such that
$x_i$ and $x_j$ are $k$-adjacent if and only if  $\vert\, i-j \,\vert=\pm1(mod\,l)$.
Then, the number $l$ of $SC_k^{n,l}$ depends on both the dimension $n$ of ${\mathbb Z}^n$ and the $k$-adjacency (see many types of $SC_k^{n,l}$ in (5) on the page of 6 of \cite{H8}).\\
For a digital image $(X,k)$ and $x \in X$, 
we follow the notation
$$N_k(x,1):=\{x^\prime \in X\,\vert \,x\,\,\text{is}\,\,k\text{-adjacent to}\,\,x^\prime\} \cup \{x\}, \eqno(2.2) $$
which is called a digital $k$-neighborhood of $x$ in $(X,k)$ \cite{H1,H2,H3,H5}.
Indeed, this notion has been effectively used in studying both pseudocovering spaces and  digital covering spaces.
For every point $x$ of a digital image $(X, k)$,
 an $N_k(x, 1)$ always exists in $(X,k)$,
the digital continuity of \cite{R1} can be represented by the following form.

\begin{proposition} \cite{H3,H5}
	Let $(X, k_0)$ and $(Y, k_1)$ be digital images in ${\mathbb Z}^{n_0}$
	and ${\mathbb Z}^{n_1}$, respectively. A~function $f: X \to Y$ is
	$(k_0, k_1)$-continuous if and only if for every point  $x\in X$,
	$f(N_{k_0}(x,1))$ is a subset of $N_{k_1}(f(x),1)$.
\end{proposition}

Owing to a digital $k$-graph theoretical feature of a digital image $(X, k)$,
we have often used a {\it $(k_0, k_1)$-isomorphism} in \cite{H4} instead of a {\it $(k_0, k_1)$-homeomorphism} in \cite{B2}, as follows:

\begin{definition} \cite{B2} (see also \cite{H4})
	For two digital images $(X, k_0)$ in ${\mathbb Z}^{n_0}$ and $(Y, k_1)$ in
	${\mathbb Z}^{n_1}$,
	a map $h : X \to Y$ is called a {\it $(k_0, k_1)$-isomorphism} if $h$
	is a $(k_0, k_1)$-continuous bijection and
	further, $h^{-1}: Y \to X$ is $(k_1, k_0)$-continuous.
	If $n_0 = n_1$ and $k_0 = k_1$, then we call it a {\it $k_0$-isomorphism}.
\end{definition}

Based on this approach, we can develop the notion of a $\lq\lq$radius $2$-$(k_0, k_1)$-isomorphism" or a radius $2$-$(k_0, k_1)$-covering map \cite{H1,H2} to establish the so-called $\lq\lq$ digital homotopy lifting theorem"
which is essential in studying digital homotopy theory.

\section{\bf Remarks on the earlier verion of a pseudo-$(k_0,k_1)$-covering space in \cite{H6,H9}}\label{s3}

Since the notions of a digital $(k_0,k_1)$-covering map and a pseudo-$(k_0,k_1)$-covering map are so related to the notion of a (weakly) local  $(k_0,k_1)$-isomorphism, we first need to recall it, as follows: 

\begin{definition} \cite{H1,H3,H9}
	For two digital images $(X, k_0)$ in ${\mathbb Z}^{n_0}$ and $(Y, k_1)$ in ${\mathbb Z}^{n_1}$, consider a  map  $h: (X,k_0) \to  (Y,k_1)$. Then~the map $h$ is said to be a local $(k_0, k_1)$-isomorphism
	if  for every $x \in X$, $h$ maps $N_{k_0}(x, 1)$ $(k_0, k_1)$-isomorphically onto $N_{k_1}(h(x), 1)$ i.e., the restriction map $h\vert_{N_{k_0}(x, 1)}:
	N_{k_0}(x, 1) \to N_{k_1}(h(x), 1))$ is a $(k_0, k_1)$-isomorphism.
	If $n_0 = n_1$ and $k_0 = k_1 $, then the map $h$ is called a local $k_0$-isomorphism.
\end{definition}

The paper \cite{H6} defined the following notion which is weaker than a local  $(k_0, k_1)$-isomorphism.

\begin{definition} \cite{H6} For two digital images $(X, k_0)$ in ${\mathbb Z}^{n_0}$ and
	$(Y, k_1)$ in ${\mathbb Z}^{n_1}$, a map $h: X \to Y$ is called a
	weakly local (WL-, for brevity) $(k_0, k_1)$-isomorphism
	if  for every $x \in X$, $h$ maps $N_{k_0}(x, 1)$ $(k_0, k_1)$-isomorphically onto $h(N_{k_0}(x, 1))\subset (Y, k_1)$, i.e., the restriction map $h\vert_{N_{k_0}(x, 1)}:
	N_{k_0}(x, 1) \to h(N_{k_0}(x, 1))$ is a $(k_0, k_1)$-isomorphism.
	In particular, if $n_0 = n_1$ and $k_0 = k_1 $, then the map $h$ is called a weakly local $k_0$-isomorphism (or a $WL$-$k_0$-isomorphism).
\end{definition}

Using this notion, the paper \cite{H6} defined  the notion of a pseudo-$(k_0, k_1)$-covering space, as follows:
\begin{definition} \cite{H6}
	Let $(E, k_0)$ and $(B, k_1)$ be digital images in  ${\mathbb Z}^{n_0}$  and  ${\mathbb Z}^{n_1}$, respectively.
	Let $p: E \to B $ be a surjection such that for any $b \in B$,\\
	(1) for some index set $M$, $p^{-1}(N_{k_1}(b, 1)) = \bigcup\limits_{i \in M}N_{k_0}(e_i, 1)$ with $e_i \in p^{-1}(b):=p^{-1}(\{b\})$;\\
	(2) if  $i, j \in M$ and $i \neq j$, then
	$N_{k_0}(e_i, 1)\cap N_{k_0}(e_j, 1)$ is an empty set; and\\
	(3) the restriction of $p$ to $N_{k_0}(e_i, 1)$ from $N_{k_0}(e_i, 1)$ to $N_{k_1}(b, 1)$ is a $WL$-$(k_0, k_1)$-isomorphism for all $i \in M$.\\
	Then the map $p$ is called a pseudo-$(k_0, k_1)$-covering map, $(E, p,
	B)$ is said to be a pseudo-$(k_0, k_1)$-covering and $(E, k_0)$ is called a pseudo-$(k_0, k_1)$-covering space over $(B, k_1)$.
\end{definition}

Based on the notion of a pseudo-$(k_0, k_1)$-covering space, 
the paper \cite{H6} referred to the map $p:({\mathbb Z}^+, 2) \to SC_k^{n,l}$ as in (1.1) for a pseudo-$(2,k)$-covering map. 
Indeed, the paper \cite{H6} made a mistake to take this map as a pseudo-$(2,k)$-covering map (see Remark 3.1 below).
By contary to the condition (1) of Definition 4, the map $p$ is not a pseudo-$(2,k)$-covering map, as follows:

\begin{remark} (Proposition 3.2 of \cite{P1})
	The map $p:({\mathbb Z}^+,2) \to SC_k^{n,l}:=(c_i)_{i \in [0, l-1]_{\mathbb Z}}, l \geq 4$, in (1.1) is not a pseudo-$(2,k)$-covering map.
\end{remark}
Since the proof of this assertion in \cite{P1} is incorrect, the paper \cite{H10}
corrected the errors.\\

To compare between a digital covering space and a pseudocovering space, we need to recall 
the notion of a digital covering space as follows:

\begin{definition} \cite{H2,H3,H5}
	Let $(E, k_0)$ and $(B, k_1)$ be digital images in  ${\mathbb Z}^{n_0}$  and  ${\mathbb Z}^{n_1}$, respectively.
	Let $p: E \to B $ be a surjection such that for any $b \in B$,
	the conditions (1) and (2) are equal to those of Definition 4; and the condition (3) is 
	the following:\\
	The restriction of $p$ to $N_{k_0}(e_i, 1)$ from $N_{k_0}(e_i, 1)$ to $N_{k_1}(p(e_i), 1)$ is a $(k_0, k_1)$-isomorphism for all $i \in M$.\\
	Then the map $p$ is called a digital $(k_0, k_1)$-covering map, $(E, p,
	B)$ is said to be a digital $(k_0, k_1)$-covering and $(E, k_0)$ is called a digital $(k_0, k_1)$-covering space over $(B, k_1)$.
\end{definition}
Based on Definitions 4 and 5, the following is obtained. 

\begin{theorem}  (see Corollary 4 of \cite{H8})
	In Definition 4, as a special case, assume that	$(E, k_0)$  and $(B, k_1)$ are 
	$k_0$- and $k_1$-connected, respectively.
	Then a digital $(k_0, k_1)$-covering map is equivalent to a local $(k_0, k_1)$-isomorphism. 
	\end{theorem}

In relation to the study between a $WL$-$(k_0, k_1)$-isomorphic surjection and a  pseudo-$(k_0, k_1)$-covering map,  there are the following two incorrect statements in \cite{H9} (see the identity (3.8) of Proposition 3.4 and Corollary 3.5 below) which were corrected in the paper \cite{H10}, as follows:

\begin{proposition} \cite{H10} (Correction of the identity of (4.2) in Proposition 4.4 of \cite{H9}) Let $p:(E, k_0) \to (B, k_1)$ be a $WL$-$(k_0, k_1)$-isomorphic surjection. Then, for any $b \in B$ with $e_i \in p^{-1}(\{b\})$, for some index set $M$ we obtain
	$$p^{-1}(N_{k_1}(b, 1)) =\bigcup\limits_{i \in M}N_{k_0}(e_i, 1) \,\,\text{with}\,\,e_i \in p^{-1}(\{b\}). \eqno(3.8)$$
	
	This statement of (3.8) was corrected as follows (see Remark 3.10 of \cite{H10}).
	$$ \bigcup\limits_{i \in M}N_{k_0}(e_i, 1) \subset p^{-1}(N_{k_1}(b, 1))  \,\,\text{with}\,\,e_i \in p^{-1}(\{b\}). \eqno(3.9)$$
\end{proposition}

\begin{corollary} (Correction of Corollary 4.5 of \cite{H9}) (1) A $WL$-local $(k_0, k_1)$-isomorphic surjection is equivalent to a pseudo-$(k_0, k_1)$-covering map of Definition 4.\\
	This statement was corrected in \cite{H10} as follows:\\
	(2) While a pseudo-$(k_0, k_1)$-covering map of Definition 4 implies a $WL$-local $(k_0, k_1)$-isomorphic surjection, 
	the converse does not hold.\\
	(3) However, with the revised version of a pseudo-$(k_0, k_1)$-covering map in Definition 4.1 of \cite{H10},
	a $WL$-local $(k_0, k_1)$-isomorphic surjection implies a new version of a pseudo-$(k_0, k_1)$-covering map  in Definition 4.1 of \cite{H10}.
\end{corollary}

\section{\bf Summary}\label{s4}

The paper \cite{H10} revised the condition (1) of the original version of a pseudo-$(k_0,k_1)$-space.
Based on this revision, it turns out that while a digital covering space implies a revised version of a pseudo-covering space in \cite{H10}, the converse does not hold.
Besides, we note that a $WL$-$(k_0,k_1)$-isomorphic surjection is equivalent to a revised version of a pseudo-$(k_0,k_1)$-map. 
Finally, since some suitable corrections on some mistakes and errors on the study of a digital covering, a pseudocovering, and a $WL$-$(k_0,k_1)$-isomorphism were made in \cite{H10}, we can find them shortly.

\end{document}